\documentclass[11pt, amsfonts]{amsart}

\usepackage{amsmath,amssymb,amsthm,color}
\usepackage[all]{xy}

\numberwithin{equation}{section}

\SelectTips{eu}{12}

\newtheorem{thm}{Theorem}[section]

\newtheorem{lem}[thm]{Lemma}

\theoremstyle{definition}

\newtheorem{eg}[thm]{Example}

\theoremstyle{remark}
\newtheorem{rem}[thm]{Remark}

\newcommand{\N}{\mathcal{N}}
\newcommand{\G}{\mathcal{G}}
\newcommand{\K}{\mathcal{K}}
\newcommand{\Q}{\mathcal{Q}}

\newcommand{\HH}{\mathcal{H}}

\newcommand{\NN}{\mathbb{N}}
\newcommand{\ZZ}{\mathbb{Z}}
\newcommand{\conn}{{\rm conn}}

\title[Neighborhood complexes and Kronecker double coverings]{Neighborhood complexes and Kronecker double coverings}

\author{Takahiro Matsushita}
\address{Department of Mathematical Sciences, University of the Ryukyus, Nishihara-cho, Okinawa 903-0213, Japan}
\email{mtst@sci.u-ryukyu.ac.jp}



\begin{document}

\baselineskip.525cm

\maketitle

\begin{abstract}
The neighborhood complex $N(G)$ is a simplicial complex assigned to a graph $G$ whose connectivity gives a lower bound for the chromatic number of $G$. We show that if the Kronecker double coverings of graphs are isomorphic, then their neighborhood complexes are isomorphic. As an application, for integers $m$ and $n$ greater than 2, we construct connected graphs $G$ and $H$ such that $N(G) \cong N(H)$ but $\chi(G) = m$ and $\chi(H) = n$. We also construct a graph $KG_{n,k}'$ such that $KG_{n,k}'$ and the Kneser graph $KG_{n,k}$ are not isomorphic but their Kronecker double coverings are isomorphic.
\end{abstract}


\section{Introduction}

The neighborhood complex was introduced by Lov\'asz in his proof of Kneser's conjecture \cite{L}. He assigned a simplicial complex $N(G)$ to a graph $G$, and showed that a certain homotopy invariant $\conn(N(G))$, called the connectivity, gives a lower bound for the chromatic number. He used this method to compute the chromatic number of the Kneser graphs $KG_{n,k}$. After that, topological methods in graph coloring problems have been studied by many authors. We refer to \cite{K} for the background of this subject.

In the study of neighborhood complexes, the following question is quite fundamental: Does the isomorphism type (homeomorphism type, or homotopy type) of $N(G)$ determine the chromatic number $\chi(G)$? Actually, this problem was negatively solved. Walker \cite{W} and Matsushita \cite{M} deal with many examples of graphs whose neighborhood complexes are homotopy equivalent but whose chromatic numbers are different. Moreover, Walker \cite{W} gave examples that for every $n \ge 2$, there are graphs $G$ and $H$ such that $\chi(G) = n$ and $\chi(H) = n+1$, but their neighborhood complexes are isomorphic.

The purpose of this paper is to improve Walker's result:

\begin{thm} \label{thm 1}
Let $m$ and $n$ be integers greater than 2. Then there are connected graphs G and H such that $\chi(G) = m$, $\chi(H) = n$, but their neighborhood complexes are isomorphic.
\end{thm}

The method employed here is different from Walker's. In this paper, we observe that the following close relation between neighborhood complexes $N(G)$ and Kronecker double coverings $K_2 \times G$ (The precise definitions will be found in Section 2).

\begin{thm} \label{thm 2}
Let $G$ and $H$ be graphs. If $K_2 \times G \cong K_2 \times H$, then $N(G) \cong N(H)$. On the other hand, if $G$ and $H$ are stiff and $N(G) \cong N(H)$, then $K_2 \times G \cong K_2 \times H$.
\end{thm}

This theorem will be proved in Section 2. Thus to prove Theorem \ref{thm 1}, it suffices to construct graphs $X(m,n)$ and $Y(m,n)$ such that $\chi(X(m,n)) = m$ and $\chi(Y(m,n)) = n$, but $K_2 \times X(m,n) \cong K_2 \times Y(m,n)$, and this will be done in Example \ref{eg 3.2}.

Theorem \ref{thm 2} asserts that the neighborhood complex is determined by its Kronecker double covering. Thus the Kronecker double covering gives a restriction on the chromatic number. In Section 3, we construct a simple graph $KG_{n,k}'$ for $n > 2k \ge 4$ such that $K_2 \times KG_{n,k}' \cong K_2 \times KG_{n,k}$ but $KG_{n,k}' \not\cong KG_{n,k}$ (Theorem \ref{thm 3.5}). By the connectivity of $N(KG_{n,k}') = N(KG_{n,k})$, we prove $\chi(KG_{n,k}') = n-2k+2$ (Theorem \ref{thm 3.6}).

Finally, we make a remark on the box complex \cite{C, MZ}. The box complex $B(G)$ is a $\ZZ/2$-space assigned to a graph, whose underlying space is homotopy equivalent to $N(G)$. Moreover, a certain $\ZZ/2$-homotopy invariant of $B(G)$, called $\ZZ/2$-index, is a lower bound for $\chi(G)$ sharper than $\conn(N(G))$ (see \cite{MZ}).

One can ask if a similar assertion to Theorem \ref{thm 1} holds for box complexes. Since $N(G) \simeq B(G)$, it is clear that $K_2 \times G \cong K_2 \times H$ implies $B(G) \simeq B(H)$. However, there are many definitions of box complexes, and these definitions are not isomorphic but only $\ZZ/2$-homotopy equivalent. Hence the isomorphism problem concerning box complexes is not so reasonable although $K_2 \times G \cong K_2 \times H$ implies $B(G) \cong B(H)$ for every definition of box complexes as far as the author knows.

On the other hand, it is meaningful to ask if $K_2 \times G \cong K_2 \times H$ implies that $B(G)$ and $B(H)$ are $\ZZ/2$-homotopy equivalent. However, the graphs constructed in Example \ref{eg 3.2} are counter examples to this question (see Remark \ref{rem 3.3}).



\section{Neighborhood complexes}

Here we review definitions and facts concerning neighborhood complexes, and show Theorem \ref{thm 2}. For a comprehensive introduction to this subject, we refer to \cite{K}.

A {\it graph} is a pair $G = (V(G) , E(G))$ consisting of a finite set $V(G)$ together with a symmetric binary relation $E(G)$ of $V(G)$. For a pair $v$ and $w$ of vertices of $G$, we write $v \sim w$ to mean $(v,w) \in E(G)$. A {\it graph homomorphism} from a graph $G$ to a graph $H$ is a map $f \colon V(G) \to V(H)$ such that $(f \times f)(E(G)) \subset E(H)$. Let $K_n$ be the graph defined by $V(K_n) = \{ 1, \cdots , n\}$ and $E(K_n) = \{ (i,j) \; | \; i \ne j\}$. The {\it chromatic number $\chi(G)$ of $G$} is the number $$\min \{ n \ge 0 \; | \; \textrm{There is a graph homomorphism $G \to K_n$}\}.$$

Let $G$ be a graph and $v$ a vertex of $G$. Let $N(v)$ be the set of vertices adjacent to $v$. The {\it neighborhood complex $N(G)$} is the simplicial complex
$$N(G) = \{ \sigma \subset V(G) \; | \; \textrm{$\sigma$ is finite and $\sigma \subset N(v)$ for some $v$}\}$$
whose underlying set is $V(G)$. Lov\'asz \cite{L} showed that if $N(G)$ is $n$-connected, then $\chi(G) > n + 2$. He used this method to determine the chromatic numbers of Kneser graphs $KG_{n,k}$ defined as follows: Let $n$ and $k$ be positive integers satisfying $n \ge 2k$. Then the {\it Kneser graph $KG_{n,k}$} is the graph defined by
$$V(KG_{n,k}) = \{ \sigma \subset \{ 1,\cdots, n\}\; | \; |\sigma| = k\}, \; E(KG_{n,k}) = \{ (\sigma, \tau) \; | \; \sigma, \tau \in V(KG_{n,k}), \; \sigma \cap \tau = \emptyset\}.$$
It is easy to see $\chi(KG_{n,k}) \le n-2k+2$. Lov\'asz showed that $N(KG_{n,k})$ is $(n-2k-1)$-connected, and hence $\chi(KG_{n,k}) = n-2k+2$.

Next we recall the definition of Kronecker double coverings. The {\it categorical product} of $G$ and $H$ is the graph $G \times H$ defined by $V(G \times H) = V(G) \times V(H)$ and $E(G \times H) = \{ ((v,w), (v',w')) \; | \; (v,v') \in E(G), (w,w') \in E(H)\}$. The {\it Kronecker double covering of $G$} is the product $K_2 \times G$. For a more detailed discussion on the Kronecker double covering, see Section 3 or \cite{IP}. The projection $K_2 \times G \to G$, $(i,v) \mapsto v$ is a covering. Here a {\it covering} means a graph homomorphism $f \colon G \to H$ such that $f|_{N(v)} \; N(v) \to N(f(v))$ is bijective for every $v \in V(G)$. It is easy to see that for a connected graph $G$, $K_2 \times G$ is connected if and only if $\chi(G) > 2$.

Now we start the proof of Theorem \ref{thm 2}. In fact, this theorem is deduced from an observation of \cite{BGZ} concerning neighborhood hypergraphs. However, we first give a direct short proof for reader's convenience. We start with the following easy observation:

\begin{lem} \label{lem 3}
$N(K_2 \times G) \cong N(G) \sqcup N(G)$
\end{lem}
\begin{proof}
For $i = 1,2$, define $f_i \colon V(G) \to V(K_2 \times G)$ by $f_i(v) = (i,v)$. Then the sum $f_1 + f_2 \colon V(G) \sqcup V(G) \to V(K_2 \times G)$ gives an isomorphism $N(G) \sqcup N(G) \to N(K_2 \times G)$.
\end{proof}

A graph $G$ is {\it stiff} if for every pair of vertices $v$ and $w$, $N(v) \subset N(w)$ implies $v = w$. Let $F(N(G))$ denote the set of facets of $N(G)$. Then the stiffness of graphs means the map $V(G) \to F(N(G))$, $v \mapsto N(v)$ is well-defined and bijective.

Before giving the proof of Theorem \ref{thm 2}, we prove the following lemma:

\begin{lem} \label{lem 4}
Let $K$ and $L$ be finite simplicial complexes. If $K \sqcup K$ and $L \sqcup L$ are isomorphic, then $K$ and $L$ are isomorphic.
\end{lem}
\begin{proof}
Let $X_1, \cdots, X_r$ be the connected components of $K$. We prove this lemma by induction on the number $r$ of connected components of $K$. The case $r = 0$ is clear.

Let $X_i'$ be a copy of $X_i$, and so $K \sqcup K = (X_1 \sqcup X'_1) \sqcup \cdots \sqcup (X_r \sqcup X'_r)$. Similarly, let $Y_1, \cdots, Y_s$ be the connected components of $L$ and so that $L \sqcup L = (Y_1 \sqcup Y'_1) \sqcup \cdots \sqcup (Y_s \sqcup Y'_s)$. Let $f \colon K \sqcup K \to L \sqcup L$ be an isomorphism. By changing indices of $Y_i$ and exchanging $Y_i$ and $Y'_i$, we can assume $f(X_1) = Y_1$. Then $f(X'_1)$ is a connected component of $L \sqcup L$ other than $Y_1$. Note that $f(X'_1)$ and $Y'_1$ are isomorphic since $f(X'_1) \cong X'_1 \cong X_1 \cong Y_1 \cong Y'_1$ Let $g \colon L \sqcup L \to L \sqcup L$ be an isomorphism which exchanges $f(X_1)$ and $Y'_1$ and fixes other components. Then we have $gf(X_1) = Y_1$ and $gf(X'_1) = Y'_1$.

Set $K' = X_2 \sqcup \cdots \sqcup X_r$ and $L' = Y_2 \sqcup \cdots \sqcup Y_s$. Then $gf$ induces an isomorphism between $K' \sqcup K'$ and $L' \sqcup L'$. By the inductive hypothesis, we have $K' \cong L'$. Since $X_1$ and $Y_1$ are isomorphic, we conclude $K = X_1 \sqcup K' \cong Y_1 \sqcup L' = L$.
\end{proof}

\noindent {\it Proof of Theorem \ref{thm 2}.} If $K_2 \times G \cong K_2 \times H$, then Lemma \ref{lem 3} implies $N(G) \sqcup N(G) \cong N(H) \sqcup N(H)$, and hence Lemma \ref{lem 4} implies $N(G) \cong N(H)$.

On the other hand, suppose $G$ and $H$ are stiff, and let $\varphi \colon V(G) \to V(H)$ be an isomorphism from $N(G)$ to $N(H)$. Define the maps $f \colon V(G) \to V(H)$ and $g \colon V(H) \to V(G)$ by $N(f(v)) = \varphi (N(v))$ and $N(g(w)) = \varphi^{-1}(N(w))$ for all $v \in V(G)$ and $w \in V(H)$. Moreover, define the maps $\tilde{f} \colon V(K_2 \times G) \to V(K_2 \times H)$ and $\tilde{g} \colon V(K_2 \times H) \to V(K_2 \times G)$ by 
$$\tilde{f}(0,v) = (0, \varphi(v)), \; \tilde{f}(1,v) = (1, f(v)), \; \tilde{g}(0,w) = (0,\varphi^{-1}(w)), \; \tilde{g}(1,w) = (1,g(w))$$
for $v \in V(G)$ and $w \in V(H)$. Then $\tilde{f}$ and $\tilde{g}$ are graph homomorphisms, and $\tilde{g}$ is the inverse of $\tilde{f}$.
\qed

\vspace{1mm}
Now we explain that Theorem \ref{thm 2} is easily deduced from an observation in \cite{BGZ} concerning neighborhood hypergraphs. To see this, we need some terminology and notation.

Recall that a {\it (multi-)hypergraph} is a pair $\HH = (V(\HH) , \HH)$ consisting of a set $V(\HH)$ together with a multi-set of $V(\HH)$, i.e. a function $\HH \colon 2^{V(\HH)} \to \NN$. The {\it neighborhood hypergraph $\N(G)$ of a graph $G$} is the multi-hypergraph on $V(G)$ whose multi-set of hyperedges is $\N(G) = \{ N(v) \; | \; v \in V(G)\}$, in other words, $\N(G) (S) = \# \{ S = N(v) \; | \; v \in V(G)\}$ for $S \in 2^{V(G)}$.


For a hypergraph $\HH$, define the bigraph representation $B_\HH$ (the precise definition of bigraphs will be found in the beginning of Section 3) as follows: the vertex set of $B_\HH$ is $V(\HH) \sqcup \HH$, and $v \in V(\HH)$ and $S \in \HH$ are adjacent if and only if $v \in S$. There is no other adjacent relation among vertices of $B_\HH$. The bigraph $B_\HH$ determines the original hypergraph $\HH$. In fact, they used this method to show that for bipartite graphs $G$ and $H$, $G \cong H$ if and only if $\N(G) \cong \N(H)$.

From the above observation of \cite{BGZ}, one can easily show Theorem \ref{thm 2} as follows: Clearly, the bigraph representation $B_{\N(G)}$ of the neighborhood hypergraph $\N(G)$ coincides with the Kronecker double covering $K_2 \times G$. This means that $K_2 \times G \cong B_{\N(G)}$ determines $\N(G)$. Since the neighborhood complex $N(G)$ is determined by $\N(G)$, we have that $K_2 \times G$ determines $N(G)$.

On the other hand, if a graph $G$ is stiff, then the neighborhood complex $N(G)$ determines the neighborhood hypergraph $\N(G)$. In fact, the multi-set of hyperedges of $\N(G)$ is the set of facets of $N(G)$ in this case. Thus if $G$ and $H$ are stiff and $N(G) \cong N(H)$, then we have $\N(G) \cong \N(H)$ and hence $K_2 \times G \cong K_2 \times H$. This completes the proof of Theorem \ref{thm 2}.

We close this section with a few remarks.

\begin{rem}
There are graphs whose neighborhood complexes are isomorphic but whose Kronecker double coverings are different. In fact, consider the 4-cycle graph $C_4$ and the path graph $P_4$ with 4 vertices. Then the neighborhood complexes of these graphs are two 1-simplices, but $K_2 \times C_4 = C_4 \sqcup C_4$ and $K_2 \times P_4 = P_4 \sqcup P_4$.
\end{rem}

\begin{rem}
Theorem \ref{thm 2} asserts that the neighborhood complex $N(G)$ is determined by the Kronecker double covering $K_2 \times G$. Thus if $N(G)$ is $n$-connected and $K_2 \times G \cong K_2 \times H$, then $N(H)$ is also $n$-connected, and hence we have $\chi(H) > n + 2$. This means that the Kronecker double covering restricts the chromatic number.

We construct graphs $KG_{n,k}'$ in Section 3 such that $K_2 \times KG_{n,k}' \cong K_2 \times KG_{n,k}$ but $KG'_{n,k} \not\cong KG_{n,k}$ for $n > 2k \ge 4$. Since $N(KG_{n,k})$ is $(n-2k-1)$-connected (see Section 2), this means $\chi(KG_{n,k}') \ge n - 2k +2$.
\end{rem}


\section{Kronecker double coverings}

In this section, we review the theory of Kronecker double coverings, and construct graphs $X(m,n)$ and $Y(m,n)$ such that $\chi(X(m,n)) =m$ and $\chi(Y(m,n)) = n$ but $K_2 \times X(m,n) \cong K_2 \times Y(m,n)$ in Example \ref{eg 3.2}. This shows Theorem \ref{thm 1}. Moreover, we construct a family of graphs $KG_{n,k}'$ such that $K_2 \times KG_{n,k} \cong K_2 \times KG_{n,k}'$ but $KG_{n,k} \not\cong KG_{n,k}'$.

We review the Kronecker double coverings from a viewpoint of bigraphs, that is, graphs with 2-colorings. For the sake of this treatment, one can obtain a simple description of the categorical equivalence given in Theorem \ref{thm 3.1}.

A {\it bigraph}\footnote{This terminology is due to \cite{BGZ}.} is a graph $X$ equipped with a 2-coloring $\varepsilon_X \colon X \to K_2$. A {\it bigraph homomorphism} is a graph homomorphism $f \colon X \to Y$ such that $\varepsilon_Y \circ f = \varepsilon_X$. Let $\G$ be the category of graphs whose morphisms are graph homomorphisms, and $\G_{/K_2}$ the category of bigraphs whose morphisms are bigraph homomorphisms. For a graph $G$, the Kronecker double covering $K_2 \times G$ is a bigraph whose 2-coloring is the 1st projection $K_2 \times G \to K_2$.

An {\it odd involution of a bigraph $X$} is a graph homomorphism (not necessarirly a bigraph homomorphism) $\tau \colon X \to X$ satisfying $\tau^2 = {\rm id}_X$ and $\varepsilon_X (\tau (v)) \ne \varepsilon_X(v)$ for every $v \in V(X)$. A typical example of odd involutions is the involution $(1,v) \leftrightarrow(2,v)$ of the Kronecker double covering $K_2 \times G$. In fact, the following theorem (Theorem \ref{thm 3.1}) asserts that every odd involution is obtained in this way.

We consider the category $\G_{/K_2}^{odd}$ defined as follows. An object of $\G_{/K_2}^{odd}$ is a pair $(X, \tau)$ consisting of a bigraph $X$ together with an odd involution $\tau$ of it. A morphism from $(X, \tau)$ to $(X' ,\tau')$ is a bigraph homomorphism $f \colon X \to X'$ which is equivariant, i.e. $\tau' \circ f = f \circ \tau$. Clearly, the Kronecker double covering gives a functor $\K \colon \G \to \G_{/K_2}^{odd}, \; G \mapsto K_2 \times G$. Moreover, we have the following theorem (see \cite{Le} for the terminology of category theory):

\begin{thm} \label{thm 3.1}
The functor $\K \colon K_2 \times (-) \colon \G \to \G_{/K_2}^{odd}$ is a categorical equivalence.
\end{thm}
\begin{proof}
We construct a quasi-inverse $\Q \colon \G_{/K_2}^{odd} \to \G$ of $\K$ as follows. For an object $(X, \tau)$ of $\G_{/K_2}^{odd}$, define the graph $X / \tau$ by $V(X / \tau) = \{ \{ x, \tau (x)\} \; | \; x \in V(X)\}$ and
$$E(X / \tau) = \{ (\alpha, \beta) \; | \; \alpha, \beta \in V(X/ \alpha), (\alpha \times \beta) \cap E(X) \ne \emptyset \}.$$
Roughly speaking, the graph $\Q(X) = X / \tau$ is the quotient of the graph $X$ by the $\ZZ / 2$-action $\tau$. Then a morphism $f \colon (X,\tau) \to (X', \tau')$ in $\G_{/K_2}^{odd}$ induces a graph homomorphism $\Q(f) \colon X / \tau \to X' / \tau'$, and hence we have a functor $\Q \colon \G_{/K_2}^{odd} \to \G$.

This functor $\Q$ is a quasi-inverse of $\K$. In fact, it is clear that $\Q \circ \K$ and $1_{\G}$ are naturally isomorphic. The natural isomorphism $1_{\G_{/K_2}^{odd}} \to \K \circ \Q$ is given by the map $f \colon X \to K_2 \times (X/ \tau)$ defined by $f(x) = (\varepsilon(x), q(x))$, where $q \colon X \to X / \tau$ is the quotient map. It is clear that $f$ is a graph isomorphism.
\end{proof}

Now we turn to the case of bipartite graphs. For a bipartite graph $X$, an involution $\tau \colon X \to X$ is {\it odd} if for every $x \in X$, there is no path with even length joining $x$ to $\tau(x)$. If $(X, \tau)$ is a bigraph with an odd involution, then $\tau$ is odd in the sense of bipartite graphs.

Let $X$ be a bipartite graph with an odd involution $\tau$. In this case, one can construct the quotient graph $X / \tau$ in the same way as the proof of Theorem \ref{thm 3.1}. Moreover, there is a 2-coloring $\varepsilon \colon X \to K_2$ such that $(X, \tau) \in \G_{/K_2}^{odd}$. Therefore by Theorem \ref{thm 3.1}, we have $K_2 \times (X /\tau) \cong X$ as graphs.

\begin{rem}
Define the category $\G'$ as follows. An object of $\G'$ is a bipartite graph $X$ together with its odd involution $\tau$. A morphism from $(X, \tau)$ to $(X',\tau')$ is a graph homomorphism $f \colon X \to X'$ satisfying $\tau' \circ f = f \circ \tau$. Then the Kronecker double covering gives a functor $\K' \colon \G \to \G'$. However, this functor is not a categorical equivalence. In fact, there is no map $f \colon G \to G$ such that $K_2 \times f = \tau$, where $\tau$ is the canonical odd involution of $K_2 \times G$.
\end{rem}

Now we are ready to prove Theorem \ref{thm 1}.

\begin{eg} \label{eg 3.2}
We construct graphs $X(m,n)$ and $Y(m,n)$ such that $K_2 \times X(m,n) \cong K_2 \times Y(m,n)$ but $\chi(X(m,n)) = m$ and $\chi(Y(m,n)) = n$. By Theorem \ref{thm 2}, this completes the proof of Theorem \ref{thm 1}.

First, set $X_1=X_2 = K_2 \times K_n$ and $Y_1 = Y_2 = K_2 \times K_m$. Define the graph $Z(m,n)$ by identifying the following vertices of $X_1 \sqcup X_2 \sqcup Y_1 \sqcup Y_2$:
\begin{itemize}
\item $(1,1) \in V(X_1)$ and $(1,1) \in V(Y_1)$.
\item $(2,1) \in V(X_1)$ and $(1,1) \in V(Y_2)$.
\item $(1,1) \in V(X_2)$ and $(2,1) \in V(Y_1)$.
\item $(2,1) \in V(X_2)$ and $(2,1) \in V(Y_2)$.
\end{itemize}
It is clear that $Z(m,n)$ is bipartite and connected. Figure 1 depicts the graph $Z(m,n)$ in the case $m = 4$ and $n = 3$.

Next we define the odd involutions $\tau_1, \tau_2$ of $Z(m,n)$. First $\tau_1$ maps $X_i$ to $X_i$ for each $i$ and $\tau_1|_{X_i}$ is the natural involution of $X_1 = X_2 =K_2 \times K_n$, flipping $K_2$. On $Y_1 \sqcup Y_2$, the involution $\tau_1$ exchanges $Y_1$ and $Y_2$, and is given by $V(Y_1) \ni (\varepsilon,x) \leftrightarrow (\varepsilon, x) \in V(Y_2)$. Similarly, $\tau_2$ maps $Y_i$ to $Y_i$ for each $i$ and $\tau_2|_{Y_i}$ is the natural involution of $K_2 \times K_m$, flipping $K_2$. On $X_1 \sqcup X_2$, the involution $\tau_2$ is given by $V(Y_1) \ni (\varepsilon,x) \leftrightarrow (\varepsilon, x) \in V(X_2)$.

Set $X(m,n) = Z(m,n)/ \tau_1$ and $Y(m,n) = Z(m,n)/ \tau_2$. To complete the proof, we need to check $\chi(X(m,n)) = m$ and $\chi(Y(m,n)) = n$. We only prove $\chi(X(m,n)) = n$ since the other is similarly shown. However, this clearly follows from the following description of $X(m,n)$: $X(m,n)$ is obtained by identifying the following vertices of $X'_1 \sqcup X'_2 \sqcup (K_2 \times K_m)$, where $X_1' = X_2' = K_m$:
\begin{itemize}
\item $1 \in V(X'_1) = V(K_m)$ and $(1,1) \in V(K_2 \times K_n)$.
\item $1 \in V(X'_2) = V(K_m)$ and $(2,1) \in V(K_2 \times K_n)$.
\end{itemize}

Figure 1 depicts the graphs $X(m,n)$ and $Y(m,n)$ in the case $m= 4$ and $n = 3$. In this figure, the involution $\tau_1$ is the reflection in the horizontal line, and the involution $\tau_2$ is the reflection in the vertical line.

\begin{figure}[t]
\begin{picture}(260,100)(0,15)
\put(20,30){\circle*{3}} \put(0,50){\circle*{3}} \put(20,50){\circle*{3}} \put(40,50){\circle*{3}}
\put(0,70){\circle*{3}} \put(20,70){\circle*{3}} \put(40,70){\circle*{3}} \put(20,90){\circle*{3}}

\put(60,10){\circle*{3}} \put(60,50){\circle*{3}} \put(80,10){\circle*{3}} \put(80,50){\circle*{3}}
\put(60,70){\circle*{3}} \put(60,110){\circle*{3}} \put(80,70){\circle*{3}} \put(80,110){\circle*{3}}

\put(120,30){\circle*{3}} \put(100,50){\circle*{3}} \put(120,50){\circle*{3}} \put(140,50){\circle*{3}}
\put(100,70){\circle*{3}} \put(120,70){\circle*{3}} \put(140,70){\circle*{3}} \put(120,90){\circle*{3}}

\put(20,30){\line(1,1){20}} \put(20,30){\line(-1,1){20}} \put(20,30){\line(0,1){20}} \put(0,50){\line(2,1){40}} \put(40,50){\line(-2,1){40}}
\put(20,50){\line(1,1){20}} \put(20,50){\line(-1,1){20}} \put(0,70){\line(1,1){20}} \put(40,70){\line(-1,1){20}} \put(20,70){\line(0,1){20}}
\put(0,50){\line(1,1){20}} \put(40,50){\line(-1,1){20}}

\put(120,30){\line(1,1){20}} \put(120,30){\line(-1,1){20}} \put(120,30){\line(0,1){20}} \put(100,50){\line(2,1){40}} \put(140,50){\line(-2,1){40}}
\put(120,50){\line(1,1){20}} \put(120,50){\line(-1,1){20}} \put(100,70){\line(1,1){20}} \put(140,70){\line(-1,1){20}} \put(120,70){\line(0,1){20}}
\put(100,50){\line(1,1){20}} \put(140,50){\line(-1,1){20}}

\put(20,30){\line(2,-1){40}} \put(20,30){\line(2,1){40}} \put(60,10){\line(1,2){20}} \put(60,50){\line(1,-2){20}} \put(80,10){\line(2,1){40}} \put(80,50){\line(2,-1){40}}

\put(20,90){\line(2,-1){40}} \put(20,90){\line(2,1){40}} \put(60,70){\line(1,2){20}} \put(60,110){\line(1,-2){20}} \put(80,70){\line(2,1){40}} \put(80,110){\line(2,-1){40}}

\put(160,55){The graph $Z(4,3)$}
\end{picture}
\end{figure}

\begin{figure}[t]
\begin{picture}(280,85)(0,-25)
\put(0,30){\circle*{3}} \put(20,10){\circle*{3}} \put(20,50){\circle*{3}} \put(40,30){\circle*{3}} \put(60,50){\circle*{3}} \put(60,10){\circle*{3}}
\put(80,10){\circle*{3}} \put(80,50){\circle*{3}} \put(100,30){\circle*{3}} \put(120,10){\circle*{3}} \put(120,50){\circle*{3}} \put(140,30){\circle*{3}}

\put(0,30){\line(1,1){20}} \put(0,30){\line(1,-1){20}} \put(0,30){\line(1,0){40}} \put(20,10){\line(0,1){40}}
\put(20,10){\line(1,1){40}} \put(20,50){\line(1,-1){40}} \put(60,10){\line(1,2){20}} \put(60,50){\line(1,-2){20}} \put(120,10){\line(0,1){40}}
\put(80,10){\line(1,1){40}} \put(80,50){\line(1,-1){40}} \put(100,30){\line(1,0){40}} \put(120,10){\line(1,1){20}} \put(120,50){\line(1,-1){20}}

\put(170,10){\circle*{3}} \put(170,50){\circle*{3}} \put(190,30){\circle*{3}} \put(210,10){\circle*{3}} \put(210,30){\circle*{3}} \put(210,50){\circle*{3}}
\put(230,10){\circle*{3}} \put(230,30){\circle*{3}} \put(230,50){\circle*{3}} \put(250,30){\circle*{3}}
\put(270,10){\circle*{3}} \put(270,50){\circle*{3}}

\put(170,10){\line(0,1){40}} \put(170,10){\line(1,1){40}} \put(170,50){\line(1,-1){40}} \put(190,30){\line(1,0){20}} \put(210,10){\line(1,1){20}}
\put(210,10){\line(1,2){20}} \put(210,30){\line(1,1){20}} \put(210,30){\line(1,-1){20}} \put(210,50){\line(1,-1){20}} \put(210,50){\line(1,-2){20}} \put(230,10){\line(1,1){40}} \put(230,50){\line(1,-1){40}} \put(270,10){\line(0,1){40}} \put(230,30){\line(1,0){20}}

\put(28,-10){The graph $X(4,3)$} \put(178,-10){The graph $Y(4,3)$}

\put(120, -28){\bf Figure 1.}
\end{picture}
\end{figure}
\end{eg}

\begin{rem} \label{rem 3.3}
The box complexes of $X(m,n)$ and $Y(m,n)$ are not $\ZZ / 2$-homotopy equivalent if $m \ne n$. To see this, we need the following fact: The box complex is a functor from the category of graphs to the category of $\ZZ / 2$-spaces, and $B(K_n)$ and $S^{n-2}$ are $\ZZ / 2$-homotopy equivalent (Proposition 5 of \cite{MZ}).

One can suppose $m < n$. Then $K_n$ is a subgraph of $Y(m,n)$ and hence there is a $\ZZ / 2$-map from $B(K_n) \simeq_{\ZZ / 2} S^{n-2}$ to $B(Y(m,n))$. If $B(X(m,n)) \simeq_{\ZZ / 2} B(Y(m,n))$, then there is a $\ZZ / 2$-map from $S^{n-2}$ to $B(X(m,n))$. However, since $\chi(X(m,n)) = m$, there is a $\ZZ / 2$-map from $B(X(m,n))$ to $B(K_m) \simeq_{\ZZ / 2} S^{m-2}$. Thus we have a $\ZZ / 2$-map from $S^{n-2}$ to $S^{m-2}$, but this contradicts the Borsuk-Ulam theorem.
\end{rem}

In the rest of this paper, we discuss a family of simple graphs $KG_{n,k}'$ which satisfies the following interesting property: The Kronecker double covering of $KG_{n,k}'$ is isomorphic to the Kronecker double covering of  $KG_{n,k}$, but $KG_{n,k}' \not\cong KG_{n,k}$ for $n > 2k \ge 4$. In the case of $n = 5$ and $k = 2$, Imrich and Pisanski \cite{IP} shows that there is a graph $G$ such that $K_2 \times G \cong K_2 \times KG_{5,2}$ but $G \not\cong KG_{5,2}$.

Let $n$ and $k$ be integers satisfying $n > 2k \ge 4$. First, let $\alpha$ be the automorphism of the $n$-point set $\{ 1, \cdots, n\}$ which exchanges $n$ and $n-1$ and fixes the remaining points. Define the odd involution $\tau$ of $K_2 \times KG_{n,k}$ by
$$(1, \sigma) \leftrightarrow (2, \alpha(\sigma))$$
for $\sigma \in V(KG_{n,k})$. Then we set $KG_{n,k}' = (K_2 \times KG_{n,k}) / \tau$.

\begin{thm} \label{thm 3.5}
$KG_{n,k}'$ is simple and $K_2 \times KG_{n,k}' \cong K_2 \times KG_{n,k}$ but $KG_{n,k}' \not\cong KG_{n,k}$.
\end{thm}
\begin{proof}
It clearly follows from Theorem \ref{thm 3.1} that $K_2 \times KG_{n,k}' \cong K_2 \times KG_{n,k}$. We show that $KG_{n,k} \not\cong KG_{n,k}'$. Since there is no vertex $x$ of $K_2 \times KG_{n,k}$ such that $x \sim \tau(x)$, $KG_{n,k}'$ is a simple graph.

First we introduce the following notation which indicates a vertex of $KG_{n,k}'$. Let $\{ i_1, \cdots, i_k\}$ be a $k$-subset of $\{ 1,\cdots, n\}$ with $i_1 < \cdots < i_k$. If $n, n - 1\not\in \{ i_1, \cdots, i_k\}$ or $\{ n-1, n\} \subset \{ i_1, \cdots, i_k\}$, we write $(i_1, \cdots, i_k)$ to indicate the vertex $\{ (1, \{ i_1, \cdots, i_k\}), (2, \{ i_1, \cdots, i_k\})\}$ of $KG'_{n,k}$. If $i_k = n-1$, then we denote by $(i_1, \cdots, i_{k-1}, \alpha)$ the vertex $\{ (1, \{ i_1, \cdots, i_k\}), (2, \alpha \{ i_1, \cdots, i_k\})\}$ of $KG_{n,k}'$, and by $(i_1, \cdots , i_{k-1}, \beta)$ the vertex $\{ (1, \alpha\{ i_1, \cdots, i_k\}), (2, \{ i_1, \cdots, i_k\})\}$ of $KG_{n,k}'$. In this notation, we have the following adjacent relation: \begin{itemize}
\item If $i_k, j_k < n -1$, then $(i_1, \cdots , i_k) \sim (j_1, \cdots , j_k)$ iff $\{ i_1, \cdots, i_k\} \cap \{ j_1, \cdots, j_k\} = \emptyset$.

\item $(i_1, \cdots, i_{k-1}, \alpha) \not\sim (j_1, \cdots, j_{k-1}, \beta)$

\item $(i_1, \cdots, i_{k-1}, \alpha) \sim (j_1, \cdots , j_{k-1}, \alpha)$ and $(i_1, \cdots , i_{k-1} , \beta) \sim (j_1 , \cdots , j_{k-1} , \beta)$ iff\\$\{ i_1, \cdots, i_{k-1}\} \cap \{ j_1, \cdots, j_{k-1}\} = \emptyset$.
\end{itemize}

Next we recall the following property of the maximum independent sets of the Kneser graphs. For $i = 1, \cdots, n$, let $A_i$ be the set of vertices of $KG_{n,k}$ which contains $i$. Recall that the Erd\H{o}s-Ko-Rado theorem \cite{EKR} states that $A_1, \cdots, A_n$ are the maximum independent sets of $KG_{n,k}$. This family of maximum independent sets of $KG_{n,k}$ clearly satisfies the following property: For a pair of $k$-subsets $\{ i_1, \cdots, i_k\}$ and $\{ j_1, \cdots, j_k\}$ of $\{ 1,\cdots, n\}$, the intersection $A_{i_1} \cap \cdots \cap A_{i_k}$ is a one point set, and if $A_{i_1} \cap \cdots \cap A_{i_k} = A_{j_1} \cap \cdots \cap A_{j_k}$, then we have $\{ i_1, \cdots, i_k\} = \{ j_1, \cdots, j_k\}$.

Now we are ready to prove $KG_{n,k}' \not\cong KG_{n,k}$. Suppose $KG_{n,k} \cong KG_{n,k}'$. For $i = 1, \cdots, n-2$, let $B_i$ be the set of vertices of $KG_{n,k}'$ containing $i$. Then each $B_i$ is a maximum independent set of $KG_{n,k}'$ since $KG_{n,k} \cong KG_{n,k}'$ and $|B_i| = \binom{n-1}{k-1}$. There are two maximum independent sets $C_1$ and $C_2$ of $KG_{n,k}'$ different from $B_1, \cdots, B_{n-2}$.

Consider the intersection $B_1 \cap \cdots \cap B_{k-1} \cap C_1$. By the above property of Kneser graphs, this determines a vertex. If $B_1 \cap \cdots \cap B_{k-1} \cap C_1 = \{ (1, \cdots, k-1, m)\}$ with $m < n -1$, then we have $B_1 \cap \cdots \cap B_{k-1} \cap B_m = B_1 \cap \cdots \cap B_{k-1} \cap C_1$, and this contradicts the above property of Kneser graphs. Hence we have $B_1 \cap \cdots \cap B_{k-1} \cap C_1 = \{ (1, \cdots, k-1, \alpha)\}$ or $\{ (1, \cdots , k-1, \beta)\}$. We assume that $B_1 \cap \cdots \cap B_{k-1} \cap C_1 = \{ ( 1, \cdots, k-1, \alpha) \}$ since the other is similarly proved. In particular, we have $(1, \cdots, k-1, \alpha) \in C_1$.

By indcution, we show $(m, \cdots, m + k -2, \alpha) \in C_1$ for $m = 1,2, \cdots, k$. Suppose that $(m, \cdots, m+k-2,\alpha) \in C_1$. Let $\{ i_1, \cdots, i_{k-1}\}$ be a $(k-1)$-subset of $\{ 1, \cdots, n-2\}$ such that $\{ m, \cdots, m + k -1\} \cap \{ i_1, \cdots, i_{k-1}\} = \emptyset$. Considering the intersection $B_{i_1} \cap \cdots \cap B_{i_{k-1}} \cap C_1$, we deduce that $(i_1, \cdots, i_k, \alpha) \in C_1$ or $(i_1, \cdots, i_k, \beta)\in C_1$ in a similar way. Since $C_1$ is independent and $(m, \cdots, m+ k -2, \alpha) \sim (i_1, \cdots, i_{k-1}, \alpha)$, we have that $(i_1, \cdots, i_{k-1}, \beta) \in C_1$. Next by considering the intersection $B_{m+1} \cap \cdots \cap B_{m+k-1} \cap C_1$, we have that $( m+1, \cdots, m+ k - 1, \alpha) \in C_1$ or $(m+1 , \cdots, m+k-1, \beta) \in C_1$. Since $C_1$ is independent and the $(i_1, \cdots, i_{k-1}, \beta) \sim (m+1, \cdots, m+k-1, \beta)$, we have  $(m + 1, \cdots, m+ k -1, \alpha) \in C_1$. Thus the induction follows.

Hence we have $(1,\cdots, k-1, \alpha), (k, \cdots, 2k-2, \alpha) \in C_1$. However, $C_1$ is independent and $(1, \cdots, k-1, \alpha) \sim (k, \cdots, 2k-2, \alpha)$. This is a contradiction.
\end{proof}

We close this paper with determining the chromatic number of $KG_{n,k}'$.

\begin{thm} \label{thm 3.6}
$\chi(KG_{n,k}') = n -2k + 2$
\end{thm}
\begin{proof}
Since $K_2 \times KG_{n,k}' \cong K_2 \times KG_{n,k}$, it follows from Theorem \ref{thm 2} that $N(KG_{n,k}') = N(KG_{n,k})$. Since $N(KG_{n,k})$ is $(n-2k-1)$-connected, we have that $\chi(KG_{n,k}') \ge n -2k+2$. So it suffices to construct an $(n-2k+2)$-coloring on $KG_{n,k}'$.

This is proved by induction on $n$. First, note that $KG_{2k,k}$ is copies of $K_2$, and hence $K_2 \times KG_{2k,k}$ is also copies of $K_2$. Since $KG_{2k,k}' = (K_2 \times KG_{2k,k}) / \tau$ is simple, we have that $KG_{2k,k}'$ is copies of $K_2$.

By the notation introduced in the proof of Theorem \ref{thm 3.5}, it is clear that $KG_{n,k}'$ is an induced subgraph of $KG_{n+1,k}'$. The set of vertices of $KG_{n+1,k}'$ not contained in $KG_{n,k}'$ is $B_{n-1}$ in the proof of Theorem \ref{thm 3.5}. Since $B_{n-1}$ is an independent set, we can construct an $(n-2k+3)$-coloring $c$ of $KG_{n+1, k}'$ as follows:
$$c(x) = \begin{cases}
c'(x) & (x \in V(KG_{n,k}')) \\
n - 2k + 3 & (x \in B_{n-1}).
\end{cases}$$
Here $c'$ is an $(n-2k+2)$-coloring of $KG_{n,k}'$.
\end{proof}

\end{document}